\documentclass[12pt]{article}
\usepackage[cp1251]{inputenc}
\usepackage[russian, english]{babel}
\usepackage{amsfonts,amssymb}
\usepackage{xcolor}
\begin{document}
\selectlanguage{russian}
\newtheorem{theor}{Теорема }
\newtheorem{defin}{Definition}
\newtheorem{cor}{Следствие}

\title{Inequalities of Jackson–Stechkin type for approximation of elements of Hilbert space (in Russian)}

\author{Vladyslav Babenko, Svitlana Konareva}

\date{Dnipropetrovsk National University}

\maketitle

\begin{center}

\leftline{\normalsize УДК 517.5}

\bigskip

{\Large Неравенства типа Джексона -- Стечкина для аппроксимации элементов гильбертова пространства}\\

\bigskip

{\Large В.Ф. Бабенко, С.В. Конарева}\\

\medskip

{\small Днепропетровский национальный университет, Днепропетровск}

\end{center}

\begin{abstract}

В роботi введенi новi характеристики елементiв гiльбертового простору -- узагальненi модулi неперервностi $\omega_\varphi(x;L_{p,V}([0,\delta]))$ i отриманi новi точнi нерiвнoстi типу Джексона -- Стечкiна з цими модулями неперервностi для апроксимацiї елементiв гiльбертового простору. Цi результати включають в себе багато вiдомих нерiвнoстей для апроксимацiї перiодичних функцiй тригонометричними полiномами, апроксимацiї неперiодичних функцiй цiлими функцiями експоненцiального типу, аналогiчнi результати для майже перiодичних функцiй та iншi. Ряд результатiв є новими вже в цих класичних ситуацiях.

\bigskip

In this paper we introduced a new characteristics of the elements of a Hilbert space - generalized moduli of continuity $\omega_\varphi(x;L_{p,V}([0,\delta]))$ and obtain new exact inequalities of Jackson - Stechkin type with these moduli of continuity for the approximation of elements of a Hilbert space. These results include many well-known inequalities for approximation of periodic functions by trigonometric polynomials, approximation of non-periodic functions by entire functions of exponential type, similar results for almost periodic functions, and others. A number of results is new even in these classic situations.

\bigskip

В роботе введены новые характеристики элементов гильбертова пространства -- обобщенные модули непрерывности $\omega_\varphi(x;L_{p,V}([0,\delta]))$ и получены новые точные неравенства типа Джексона -- Стечкина с этими модулями непрерывности для аппроксимации элементов гильбертова пространства. Эти результаты включают в себя много известных неравенств для аппроксимации периодических функций тригонометрическими полиномами, аппроксимации непериодических функций целыми функциями експоненциального типа, аналогичные результаты для почти периодических функций и другие. Ряд результатов являются новыми уже в этих классических ситуациях.

\end{abstract}

\newpage

{
{1. \bf Введение.}
{Пусть $X$ -- нормированное пространство над полем комплексных чисел.
Наилучшим приближением элемента $x \in X$
подпространством $W \subset X$ называется величина}
$$
E \left(f, W\right)_X = \inf\limits_{h \in W} \left\|f -
h\right\|_X.
$$
Если $G$ есть действительная ось $\mathbb{R}$  или единичная
окружность $\mathbb{T}$, реализо\-ванная как отрезок $\left[0, 2\pi
\right]$ с отождествлёнными концами, то через $X(G)$ будем обозначать нормированное пространство комплекснозначных функций, заданных на $G$. В частности, мы будем иметь дело с пространствами $L_p(G)$ и $C(G)$.
Наилучшее приближение функции $f\in X(\mathbb{T})$ тригонометрическими полиномами
по\-рядка не выше $n-1$ будем обозначать через
$E_{n}(f)_{X(\mathbb{T})}$.

Модуль непрерывности порядка $m\in \mathbb{N}$
функции $f\in X(G)$ определяется так:
\[
\omega_m \left(f, \delta\right)_{X(G)} = \sup\limits_{0\leq t
\leq \delta}\left\|\Delta_{t}^{m} f\right\|_{X(G)}=\sup\limits_{0\leq t
\leq \delta}\left\|\sum\limits_{k = 0}^{m}(-1)^{m-k} C_{m}^{k}
f(\cdot + kt)\right\|_{X(G)}\qquad .
\]

Неравенства, оценивающие величины $E \left(f, W\right)_{X(G)}$ через
значение модуля непрерывности $\omega_m \left(f, \delta\right)_{X(G)}$ в
некоторой точке $\delta$, называются неравенcтвами типа Джексона
(Джексона -- Стечкина при $m\geq 2$).
Первое точное неравенство типа Джексона для наилучших равномерных приближений функций из
$C(\mathbb{T})$ тригонометрическими полиномами было получено Н.П. Корнейчуком~\cite{Korn} в 1962 году.
Аналогичный ре\-зультат для наилучших равномерных
приближений функций $f \in C(\mathbb{R})$ целы\-ми функциями
экспоненциального типа $\sigma$ был получен В.К. Дзядыком~\cite{Dz}.

В 1967г. Н.И. Черных~\cite{Cher, Cher N} доказал следующие два
неулучшаемых нера\-венства для функций
$f\in L_2\left(\mathbb{T}\right)$:
\begin{equation}\label{2}
 E_{n}\left(f\right)_{L_2 \left(\mathbb{T}\right)}\leq
\frac{1}{\sqrt{2}} \omega_1\left(f, \frac{\pi}{n}\right)_{L_2
\left(\mathbb{T}\right)},\;\;\;E_{n}\left(f\right)_{L_2 \left(\mathbb{T}\right)}\leq
\frac{1}{\sqrt{C_{2m}^m}} \omega_m\left(f,
\frac{2\pi}{n}\right)_{L_2 \left(\mathbb{T}\right)},  m\geq 2.
\end{equation}
Аналогичные результаты для наилучших $L_2\left(\mathbb{R}\right)$-приближений функций $f \in L_2\left(\mathbb{R}\right)$ целыми функциями экспоненциального типа $\sigma$ были получены в~\cite{Ibr, Popov}, а для наилучших приближений
$B^{2}$-почти периодических функций в
\cite{Pritula,PritulaYA}.

Для получения неравенств (\ref{2}) в~\cite{Cher, Cher N} были установлены точные неравенства вида
\begin{equation}\label{3}
  E_n(f)\le K\left\{\int\limits_0^\delta\omega^2_m(f,t)v(t)dt\right\}^{1/2}
\end{equation}
c $\delta=\pi/n$ и $v(t)=\sin nt$ для $m=1$ и с $\delta=2\pi/n$ и $v(t)=\sin nt+2\sin(nt/2)$ для $m\ge 2$.

 Л. В. Тайков~\cite{Taikov1, Taikov2} начал систематическое исследование задачи о точных неравенствах вида (\ref{3}). В дальнейшем работы многих математиков были посвящены получению точных неравенств такого типа. Информацию о полученных в этом направлении результатах и дальнейшие ссылки см. в работах  \cite{Shabozov, Vakarchuk}.

В связи сo вторым из неравенств (\ref{2}) изучались, в частности, следующие вопросы.

1. Чему равна точная константа в неравенстве
\begin{equation}\label{4}
 E_{n}\left(f\right)_{L_2 \left(\mathbb{T}\right)}\leq
\chi \omega_m\left(f, \frac{\pi}{n}\right)_{L_2
\left(\mathbb{T}\right)}, m\geq 2?
\end{equation}

2. Каково минимальное $\delta > 0$ такое, что для произвольной
функции $f \in L_2 \left(\mathbb{T}\right)$
\begin{equation}\label{5}
 E_{n}\left(f\right)_{L_2 \left(\mathbb{T}\right)}\leq
\frac{1}{\sqrt{C_{2m}^m}} \omega_m\left(f,
\frac{\delta}{n}\right)_{L_2 \left(\mathbb{T}\right)},  m\geq 2?
\end{equation}

В~\cite{Vas, Step} показано, что для точной константы $\chi$ в
неравенстве (\ref{4}) справедлива оценка
$
\chi \leq \frac{\sqrt{m+1}}{2^m}.
$
В~\cite{Vas} также показано, что для
минимальной точки $\delta$ в неравенстве (\ref{5}) справедлива
оценка
$
\delta \leq {1.4 \pi}.
$

Неравенства типа Джексона-Стечкина с обобщенными
модулями непрерывности, введенными в работах Шапиро и Бомана \cite{Shap, Boman, Bom} (совокупность таких модулей непрерывности включает в себя, помимо классических, модули непрерывности,} по\-рожденные более
общими, по сравнению с $\Delta_t^m f$ конечно разностными
опера\-торами, разностными
операторами дробных порядков и многие другие), изучались, например, в работах ~\cite{Babenko, Vas, Vasil, Kozko A}.

В ряде работ (см., например,~\cite{Gorb,BabSav2}) задача о
точных неравенствах типа Джексона-Стечкина изучалась в
абстрактных гильбертовых пространствах. Целью данной статьи является дальнейшее изучение этой задачи.

 В \S 2 мы приводим два неравенства для операторов в гильбертовом пространстве. Некоторые необходимые факты из спектральной теории приведены в \S 3. В \S 4 вводятся обобщенные модули непрерывности элементов гильбертова пространства.  Точные оценки величин $|(x-\Lambda x, f)|$ в терминах введенных характеристик получены в \S 5 (при этом используются неравенства из \S 2), а неравенства типа Джексона -- Стечкина в гильбертовом пространстве -- в \S 6. Некоторые конкретизации результатов \S 6 обсуждаются в \S 7.

{2. \bf Одно неравенство для операторов в гильбертовом пространстве.} Пусть $H$ - комплексное гильбертово пространство со скалярным
произведением $(\cdot, \cdot)$ и нормой $\|\cdot\|= \sqrt{(\cdot,
\cdot)}$. Пусть $S,T:H\to H$ -- линейные ограниченные операторы, такие, что $ST=TS$. Будем предполагать, что оператор $S|_{T(H)}\; :\; T(H)\to S(T(H))$ имеет обратный $(S|_{T(H)})^{-1}$. Через $T^*$ и $S^*$ будем, как обычно, обозначать сопряженные операторы.

\begin{theor}
  Для любых $f\in H$ и $x\in H$ справедливо неравенство
  \begin{equation}\label{6}
    |(Tx,f)|\le \| Sx\|\cdot\left\| \left((S|_{T(H)})^{-1}T\right)^*f\right\|.
  \end{equation}
 Если
 \begin{equation}\label{condition}
 \left((S|_{T(H)})^{-1}T\right)^*(H)\subset S(T(H)),
 \end{equation}
то неравенство (\ref{6}) является точным и обращается в равенство для
  \begin{equation}\label{7}
  \tilde{x}=(S|_{T(H)})^{-1}\left((S|_{T(H)})^{-1}T\right)^*f.
  \end{equation}
\end{theor}

$\Box$
Имеем
\[
|(Tx,f)|=|((S_{T(H)})^{-1}S|_{T(H)}Tx,f)|=|((S_{T(H)})^{-1}STx,f)|=
\]
\[
=|((S_{T(H)})^{-1}TSx,f)|=\left|\left(Sx, \left((S|_{T(H)})^{-1}T\right)^*f\right)\right|\le \| Sx\|\cdot\left\| \left((S|_{T(H)})^{-1}T\right)^*f\right\|.
\]
Неравенство (\ref{6}) доказано.

Для элемента $\tilde{x}$ имеем
\[
S\tilde{x}=S(S|_{T(H)})^{-1}\left((S|_{T(H)})^{-1}T\right)^*f=S|_{T(H)}(S|_{T(H)})^{-1}\left((S|_{T(H)})^{-1}T\right)^*f=\left((S|_{T(H)})^{-1}T\right)^*f.
\]
Следовательно,
\[
(T\tilde{x},f)=\left(S\widetilde{x}, \left((S|_{T(H)})^{-1}T\right)^*f\right)=\left(\left((S|_{T(H)})^{-1}T\right)^*f, \left((S|_{T(H)})^{-1}T\right)^*f\right)=
\]
\[
=\left\|\left((S|_{T(H)})^{-1}T\right)^*f\right\|^2=\left\|\left((S|_{T(H)})^{-1}T\right)^*f\right\|\cdot \left\|\left((S|_{T(H)})^{-1}T\right)^*f\right\|=
\]
\[
=\| S\tilde{x}\|\cdot \left\|\left((S|_{T(H)})^{-1}T\right)^*f\right\|.
\]
Теорема доказана.
$\Box$

\begin{cor}
Для любого $x\in H$
\begin{equation}\label{ner_dlya_norm}
\|Tx\|\le \| Sx\|\cdot\left\| \left((S|_{T(H)})^{-1}T\right)^*\right\|.
\end{equation}
Если выполнено условие (\ref{condition}), то неравенство (\ref{ner_dlya_norm}) является точным.
\end{cor}

$\Box$ Из неравенства (\ref{6}) получаем
\[
\|Tx\|=\sup\limits_{f\in H\atop \| f\|\le 1}|(Tx,f)|\le\| Sx\|\sup\limits_{f\in H\atop \| f\|\le 1}\left\|\left((S|_{T(H)})^{-1}T\right)^*f\right\|=\| Sx\|\cdot\left\| \left((S|_{T(H)})^{-1}T\right)^*\right\|.
\]
Докажем точность полученного неравенства. Пусть элемент $\tilde{f}\in H$, $\| \tilde{f}\|=1$, таков, что
\[
\left\| \left((S|_{T(H)})^{-1}T\right)^*\tilde{f}\right\|=\left\| \left((S|_{T(H)})^{-1}T\right)^*\right\|.
\]
Положим
\[
\tilde{x}=(S|_{T(H)})^{-1}\left((S|_{T(H)})^{-1}T\right)^*\tilde{f}.
\]
Тогда будем иметь (см. доказательство точности неравенства (\ref{6}))
\[
\| T\tilde{x}\|\ge |(T\tilde{x},\tilde{f})|=\| S\tilde{x}\|\cdot\|\left((S|_{T(H)})^{-1}T\right)^*\tilde{f}\|=\| S\tilde{x}\|\cdot\| \left((S|_{T(H)})^{-1}T\right)^*\|.
\]
Следствие доказано.
$\Box$

\bigskip

{3. \bf Необходимые сведения из спектральной теории операторов в гильбертовом пространстве.}
Пусть задано гильбертово пространство $H$. Говорят (см.~\cite[гл. ХIII, \S 1]{Berezanskiy}) что на $\sigma$-алгебре ${\cal B}$ борелевских подмножеств числовой прямой задано разложение единицы $E$, если каждому $\beta\in {\cal B}$ поставлен в соответствие проектирующий оператор $E(\beta)$ в $H$ причем выполняются следующие условия:
\begin{enumerate}
  \item $E(\emptyset)=0,\;\; E(\mathbb{R})=I$;
  \item Для любой последовательности $\{\beta_j\}_{j=1}^\infty\subset {\cal B}$, состоящей из попарно непересекающихся множеств,
      \[
      E\left(\bigcup\limits_{j=1}^\infty\beta_j\right)=\sum\limits_{j=1}^\infty E(\beta_j).
      \]
\end{enumerate}

Заданное разложение единицы порождает ~\cite[гл. ХIII, \S\S 6, 7]{Berezanskiy} группу унитарных операторов $U_t$ и самосопряженный оператор $A$:

$$
U_t x = \int\limits_{- \infty}^{+ \infty} e^{ist} dE(s) x, \;\;\;
(t \in \mathbb{R}),\;\;\;\; Ax=\int\limits_{- \infty}^{+ \infty} t dE(t) x.
$$
Измеримая и почти везде конечная функция $F:\mathbb{R}\to \mathbb{C}$ определяет функцию $F(A)$ от самосопряженного оператора $A$:
\[
F(A)x=\int\limits_{-\infty}^\infty F(t)dE(t)x.
\]
При этом
\[
D(F(A)):=\left\{ x\in H\; :\; \int\limits_{-\infty}^\infty |F(t)|^2d(E(t)x,x)<\infty\right\}\;\;\mbox{и}\;\;\| F(A)\|^2=\int\limits_{-\infty}^\infty |F(t)|^2d(E(t)x,x).
\]
Для сопряженного оператора $F(A)^*$ и для оператора, обратного к $F(A)$ (если он определен) имеем
\[
F(A)^*x=\int\limits_{-\infty}^\infty \overline{F(t)}dE(t)x\qquad\mbox{и}\qquad F(A)^{-1}x=\int\limits_{-\infty}^\infty \frac 1{F(t)}dE(t)x.
\]

{\bf 4. Обобщенные модули непрерывности элементов гильбертова пространства.} Обозначим через $\Phi$ множество непрерывных неотрицательных $2\pi$ - периодических функций $\psi$, имеющих нигде не плотное множество нулей и таких, что $\psi (0)=0$. Пусть
$\widehat{f_s}$ - коэффициенты Фурье функции $f$,
$\psi(\cdot)\in \Phi$.
В работах  ~\cite{Shap, Boman, Bom} было предложено
обобщенным модулем непрерывности функции $f \in L_2
\left(\mathbb{T}\right)$  называть функцию
$$
\omega_\psi\left(f, \delta\right)_{L_2(\mathbb{T})} = \max\limits_{t \in
[0,\delta]}\left(\sum\limits_{s \in
\mathbb{Z}}\psi(st)|\widehat{f_s}|^2\right)^{\frac{1}{2}},\qquad \delta\ge 0.
$$

Пусть $\varphi : \mathbb{C}\to \mathbb{C}$ -- непрерывная функция, такая, что $\psi(t)=|\varphi(e^{it})|^2\in \Phi$. В частности, $\varphi(1)=0$ и на любой дуге окружности $|z|=1$ функция $\varphi(z)$ отлична от тождественного нуля.
Определим обoбщенную разность элемента $x\in H $ с шагом $t$ полагая
$$
\Delta_{t}^{\varphi}x= \int\limits_{- \infty}^{+
\infty}\varphi\left(e^{its}\right)dE(s)x.
$$
Ясно, что
$$
\left\|\Delta_{t}^{\varphi}x\right\|^2 = \int\limits_{- \infty}^{+
\infty}\left|\varphi\left(e^{its}\right)\right|^2d\left(E(s)x,x\right)
= \int\limits_{- \infty}^{+
\infty} \psi\left(st\right)d\left(E(s)x,x\right).
$$
Отметим, что $\left\|\Delta_{t}^{\varphi}x\right\|$ непрерывно зависит от $t$ и $\left\|\Delta_{t}^{\varphi}x\right\|\to 0, \; t\to 0$.

Обобщенным модулем непрерывности элемента $x$ гильбертова
пространства $H$ назовем
\begin{equation}\label{8}
\omega_{\varphi}\left(x, \delta\right)
=\max\limits_{0\leq t\leq\delta}\left\|\Delta_{t}^{\varphi}x\right\|=\left\|\left\|\Delta_{t}^{\varphi}x\right\|\right\|_{C([0,\delta])}
=\left\|\int\limits_{-\infty}^{+\infty}\psi\left(st\right)d\left(E(s)x,x\right)\right\|_{C([0,\delta])}^{1/2}.
\end{equation}

 Кроме $\omega_{\varphi}\left(x, \delta\right)=\omega_{\varphi}\left(x, C([0,\delta])\right)$, будем рассматривать характеристики $\omega_\varphi(x;L_{p,V}([0,\delta]))$, в которых $1\le p<\infty$, $V(t)$ есть вес, то есть неотрицательная интегрируемая на $[0,1]$ функция, отличная от нуля на множестве полной меры.
Положим
\[
\omega_\varphi(x;L_{p,V}([0,\delta]))=\left(\frac 1\delta\int\limits_0^\delta \|\Delta^\varphi_t x\|^pV\left(\frac t\delta\right)dt\right)^{1/p}.
\]
Ясно, что
\[
\omega_\varphi(x;L_{2,V}([0,\delta]))=\left(\frac 1\delta\int\limits_0^\delta V\left( \frac t\delta\right)\int\limits_{- \infty}^{+
\infty}\psi\left(st\right)d\left(E(s)x,x\right)dt\right)^{1/2}=\left(\int\limits_{- \infty}^{+
\infty}\Gamma(V,\delta t)d\left(E(s)x,x\right)\right)^{1/2},
\]
где
$$
\Gamma (V;t)=\int\limits_0^1 \psi(ts)V(s)ds.
$$
Функция $\Gamma(V;t)$ непрерывно зависит от $t$ и обращается в нуль только в точке нуль.
Отметим, что при любом $\delta >0$ будет $\omega_\varphi(x;L_{p,V}([0,\delta]))\to \omega_\varphi(x;C([0,\delta]))$, если $p\to \infty$.

Ниже  нам будет удобно предполагать, что $\|V(\cdot)\|_1:=\| V(\cdot)\|_{L_1([0,1])}=1$. Тогда, в силу неравенства Гельдера, для $2\le p\le\infty$ при всех $\delta >0$
 \begin{equation}
  \label{gelder}
 \omega_\varphi(x;L_{2,V}([0,\delta]))\le \omega_\varphi(x;L_{p,V}([0,\delta])).
 \end{equation}

\bigskip

{\bf 5. Оценки величин $|(x-\Lambda x, f)|$.}
Будем рассматривать задачи аппроксимации элементов гильбертова пространства подпространствами вида
$$
W_{\sigma}= \left\{\int\limits_{|t|<\sigma} dE(s) g \;\;\; :
\;\;\; g\in H\right\}, \;\;\; \sigma > 0.
$$

Для аппроксимации будем использовать линейные методы приближения вида
$$
\Lambda x = \int\limits_{|t|<\sigma} \lambda (t) dE(t) x,
$$
где $\lambda (t)$ - непрерывная в $(-\sigma,\sigma)$, ограниченная, комплекснозначная функция, тождественно равная единице в некотором интервале $(-\epsilon,\epsilon),\; 0<\epsilon<\sigma$. Тогда
$$
(I-\Lambda)x=x - \Lambda x = \int\limits_{|t|< \sigma}\left(1- \lambda (t)
\right)dE(t) x +  \int\limits_{|t|\geq \sigma} dE(t)
x=\int\limits_{-\infty}^\infty\theta (t)dE(t)x,
$$
где $\theta (t)=1-\lambda (t)$, если $|t|<\sigma$, и $\theta (t)=1$,
если $|t|\ge\sigma$.


Для любого элемента $x\in H$ такого, что $x\neq U_tx$ при некотором
$t$, рассмотрим значение функционала $f\in H^*=H$ на разности
$x-\Lambda x $. Получим неравенства, связывающие $|(x-\Lambda x,f)|$ и $\omega_\varphi(x;L_{2,V}([0,\delta])).$
В теореме 1 положим $S=\Gamma(V,\delta \cdot)^{1/2}(A)$, $T=I-\Lambda$. $T$ и $S$ -- ограниченные операторы. При этом, в силу того, что $\theta(t)\equiv 0$ для $t\in (-\epsilon,\epsilon)$, а непрерывная функция  $\Gamma(V,\delta t)$ обращается в нуль только в точке нуль, оператор $\left( S_{T(H)}\right)^{-1}$ существует,
\[
 S|_{T(H)}^{-1}T=\int\limits_{-\infty}^\infty \frac{\theta(t)dE(t)}{\Gamma(V,\delta t)^{1/2}},\;\;\;\left(\left(S|_{T(H)}\right)^{-1}T\right)^*=\int\limits_{-\infty}^\infty \frac{\overline{\theta(t)}dE(t)}{{\Gamma(V,\delta t)^{1/2}}}
\]
(как обычно, считаем, что $0/0=0$) и для любых $x,f\in H$
\begin{equation}\label{22}
\| Sx\|^2=\omega_\varphi(F(A)x;L_{2,V}([0,\delta]))^2,\;\;\;
\left\| \left(\left( S_{T(H)}\right)^{-1}T\right)^*f\right\|^2=\int\limits_{-\infty}^\infty \frac{|\theta(t)|^2d(E(t)f,f)}{\Gamma(V,\delta t)}.
\end{equation}
Используя неравенство (\ref{6}), получим
\begin{equation}
 \label{9} |(x-\Lambda x, f)|\leq  \left(\;\int\limits_{-\infty}^\infty\frac {|\theta (t)|^2
d(E(t)f,f)}{\Gamma (V;\delta t)}\right)^{\frac 12} \omega_\varphi(x;L_{2,V}([0,\delta])).
\end{equation}
Как следует из  теоремы 1
это неравенство обращается в равенство для элемента
$$
\tilde{x}=\int\limits_{-\infty}^\infty\frac {\overline{\theta (t)}
}{\Gamma (V;\delta t)}dE(t)f.
$$
Итак, нами доказана

\begin{theor} Для произвольной непрерывной в $(-\sigma,\sigma)$, ограниченной,
комплекснозначной функции $\lambda (t)$ такой, что $\lambda(t)\equiv 1$ в $(-\epsilon,\epsilon),\; 0<\epsilon<\sigma$ , линейного метода
приближения $\Lambda x = \int\limits_{|t|<\sigma}\lambda(t)
dE(t) x,$ любого элемента $x\in H$ такого, что $x\neq U_tx$ для
некоторого $t$, любого элемента $f\in H$ и любого  веса $V(t)$ имеет место неулучшаемое
неравенство (\ref{9}).
В частности,
$$
\left|\left(x-\int\limits_{|t|<\sigma} dE(t)x,
f\right)\right|\leq \left(\;\int\limits_{|t|\geq \sigma}\frac
{d(E(t)f,f)}{\Gamma (V;\delta t)}\right)^{1/2} \omega_\varphi(x;L_{2,V}([0,\delta])).
$$
\end{theor}

{\bf 6. Неравенства типа Джексона -- Стечкина.}
Из неравенства (\ref{ner_dlya_norm}) с учетом (\ref{22}) выводим
$$
\|x-\Lambda x\|^2  \leq\left\| \left(\left( S_{T(H)}\right)^{-1}T\right)^*\right\|^2\cdot \omega_\varphi(x;L_{2,V}([0,\delta]))^2.
$$
Положим
$$
{\cal H} (V,\delta,\sigma)=
\inf\limits_{|t|\leq \sigma}\Gamma (V;\delta t),\;\;\; {\cal G}(V,\delta,\sigma)=
\inf\limits_{|t|\geq \sigma}\Gamma (V;\delta t).
$$
Будем иметь
$$
\left\| \left(\left( S_{T(H)}\right)^{-1}T\right)^*\right\|^2=\sup\limits_{\|f\|=1}\int\limits_{-\infty}^\infty\frac
{|\theta (t)|^2 d(E(t)f,f)}{\Gamma (V;\delta t)}\le
$$
$$
\le\sup\limits_{\|f\|=1}\left(\sup\limits_{|t|< \sigma}\frac
{|1-\lambda(t)|^2 }{\Gamma (V;\delta t)}\int\limits_{|t|< \sigma}d(E(t)f,f)+
\sup\limits_{|t|\geq \sigma}\frac {1}{\Gamma (V;\delta t)}\int\limits_{|t|\geq \sigma}d(E(t)f,f)\right)\le
$$
$$
\le\max\left\{\frac {|1-\lambda(t)|^2
}{{\cal H}(V;\delta; \sigma)}, \frac {1}{{\cal G}(V;\delta;\sigma)}
\right\}.
$$

Покажем, что если разложение единицы таково, что $E([t,t+\varepsilon])
\neq 0$ для произвольных $t\in \mathbb{R}$ и $\varepsilon > 0$, то
\begin{equation}\label{norma}
\left\| \left(\left( S_{T(H)}\right)^{-1}T\right)^*\right\|^2=\max\left\{\frac {|1-\lambda(t)|^2
}{{\cal H}(V;\delta; \sigma)}, \frac {1}{{\cal G}(V;\delta;\sigma)}
\right\}.
\end{equation}

Рассмотрим случай, когда
$
\max\left\{\frac {|1-\lambda(t)|^2
}{{\cal H}(V;\delta; \sigma)}, \frac {1}{{\cal G}(V;\delta;\sigma)}
\right\}=\frac {1}{{\cal G}(V;\delta;\sigma)}.
$
Зададим произвольное $\varepsilon>0$. Пусть $t_\varepsilon ,\; |t_\varepsilon |\ge \sigma ,$ таково, что
$
\Gamma (V,\delta t_\varepsilon )\le {\cal G}(V;\delta;\sigma)+\varepsilon
$
(для определенности будем считать, что $t_\varepsilon \ge \sigma .$) Пусть $\gamma >0$ настолько мало, что в интервале $[t_\varepsilon,t_\varepsilon +\gamma]$ выполняется неравенство $\Gamma (V,\delta t)\le\Gamma (V,\delta t_\varepsilon )+\varepsilon$ (в силу непрерывности функции $\Gamma (V,\delta t)$ такой выбор $\gamma$ возможен). Выберем элемент $f_\epsilon\in E([t_\varepsilon,t_\varepsilon +\gamma])(H)$ такой, что $\| f_\varepsilon\|=1$.  Будем иметь
\[
\sup\limits_{\|f\|=1}\int\limits_{-\infty}^\infty\frac
{|\theta (t)|^2 d(E(t)f,f)}{\Gamma (V;\delta t)}\ge \int\limits_{[t_\varepsilon,t_\epsilon +\gamma]}\frac
{ d(E(t)f_\varepsilon,f_\epsilon)}{\Gamma (V;\delta t)}\ge
\]
\[
\ge \frac 1{{\cal G}(V;\delta;\sigma)+2\varepsilon}\int\limits_{[t_\varepsilon,t_\epsilon +\gamma]} d(E(t)f_\varepsilon,f_\epsilon)=\frac 1{{\cal G}(V;\delta;\sigma)+2\varepsilon}.
\]
В силу произвольности $\varepsilon$ в рассматриваемом случае соотношение (\ref{norma}) доказано. Cлучай, когда $
\max\left\{\frac {|1-\lambda(t)|^2
}{{\cal H}(V;\delta; \sigma)}, \frac {1}{{\cal G}(V;\delta;\sigma)}
\right\}=\frac {1}{{\cal H}(V;\delta;\sigma)}
$,
разбирается аналогично.

Таким образом, доказана
\begin{theor} {Для любого элемента $x\in H$ такого, что
$x\neq U_tx$ при некотором $t$, справедливо неравенство}
\begin{equation}\label{teor_3_1}
\|x-\Lambda x\|^2  \leq \max\left\{\frac
{|1-\lambda(t)|^2 }{{\cal H} (V;\delta;\sigma)}, \frac
{1}{{\cal G} (V,\delta,\sigma)} \right\}\omega_\varphi(x;L_{2,V}([0,\delta]))^2.
\end{equation}
В частности, для наилучшего приближения элемента $x \in H$
подпространством $W_\sigma$ имеем
\begin{equation}\label{teor_3_2}
E_{\sigma}\left(x\right)^2=\left\| x-\int\limits_{|t|<\sigma} dE(t)x\right\|^2 \le
\frac 1{{\cal G}(V,\delta,\sigma)}\omega_\varphi(x;L_{2,V}([0,\delta]))^2.
\end{equation}
{Если разложение единицы таково, что $E([t,t+\varepsilon])
\neq 0$ для произвольных $t\in \mathbb{R}$ и $\varepsilon > 0$, то
неравенства (\ref{teor_3_1}) и (\ref{teor_3_2}) являются точными.}
\end{theor}

\begin{cor}
 В условиях теоремы 3 при $2\le p\le\infty$ имеет место неравенство
\begin{equation}\label{est}
E_{\sigma}\left(x\right) \leq
\frac{1}{{\cal G}(V,\delta,\sigma)^{\frac 12}}\omega_\varphi(x;L_{p,V}([0,\delta])).
\end{equation}
\end{cor}

 Приведем обобщения некоторых результатов из \cite{Vas,Vasil}.

Если найдется такой вес $V=\widetilde{V}$ (напомним, что $\| {V}\|_1=1$), что для него при некотором $\gamma > 0$
выполняется неравенство
\begin{equation}\label{below}
{\cal G}({V},\frac \gamma\sigma,\sigma) ={\cal G}({V},1,\gamma)\geq
{\cal I}\left(\psi\right),
\end{equation}
где ${\cal I}(\psi)=\frac 1{2\pi}\int\limits_0^{2\pi}\psi(t)dt$, то из (\ref{est}) при $2\le p\le \infty$ будем иметь
\begin{equation}\label{vas0}
 E_{\sigma}\left(x\right)\leq \frac{1}{\sqrt{{\cal I}\left(\psi\right)}}
  \omega_\varphi\left(x;L_{p,\widetilde{V}}\left(\left[0,\frac{\gamma}{\sigma}\right]\right)\right).
\end{equation}
В частности, при $p=\infty$
\begin{equation}\label{vas}
 E_{\sigma}\left(x\right)\leq \frac{1}{\sqrt{{\cal I}\left(\psi\right)}}
  \omega_\varphi\left(x;C\left(\left[0,\frac{\gamma}{\sigma}\right]\right)\right)= \frac{1}{\sqrt{{\cal I}\left(\psi\right)}}
  \omega_{\varphi}\left(x; \frac{\gamma}{\sigma}\right).
\end{equation}
В ~\cite{Vasil} доказано, что такой вес
существует, и приведена схема построения функции $\widetilde{V}(t)$. Таким образом установлено следующее утверждение.

\begin{theor}
Для любой функции $\varphi :\mathbb{C}\to \mathbb{C}$ такой, что $\psi\in \Phi$, существует точка $\gamma>0$ такая, что для любого $x\in H$ и любого $\sigma>0$ выполняются неравенства (\ref{vas0}) и (\ref{vas}).
\end{theor}

Сузим по сравнению с $\Phi$ класс рассматриваемых функций $\psi$. Через $\Psi$ обозначим (см. \cite{Vas}) совокупность функций $\psi\in \Phi$ таких, что

\begin{enumerate}
  \item $\psi(-t)=\psi(t)$ и $\psi(\pi -t)=\psi(\pi +t)$ для $t\in\mathbb{R}$,
  \item $\frac 1t\int\limits_0^t\psi(s)ds\le \frac 1\pi\int\limits_0^\pi\psi(s)ds$ для любого $t\in (0,\pi)$.
\end{enumerate}

Пусть
\[
Z(t)=\left\{ \begin{array}{clcr}\mbox{$\frac{2 t}{7}$}, &\mbox{$t\in\left[0, \frac{1}{7}\right]$},\\
\mbox{$\frac{-t^2}{2}+\frac{3 t}{7}-\frac{1}{98}$},&\mbox{$t\in\left[\frac{1}{7}, \frac{5}{7}\right]$},\\
\mbox{$\frac{t^2}{2}-t+\frac{1}{2}$},&\mbox{$t\in\left[\frac{5}{7},1\right].$}
\end{array}\right.
\]
Положим $V^*(t)=Z(t)/\| Z(\cdot)\|_1$. Из результатов работы ~\cite[теорема 1]{Vas} следует, что для веса $V=V^*$, для $\psi\in \Psi$ и любого $\gamma\ge \frac{7\pi}{5}$ выполняется неравенство (\ref{below}). Поэтому справедлива

\begin{theor}
 Для любой функции $\varphi :\mathbb{C}\to \mathbb{C}$ такой, что $\psi\in \Psi$, для любого $\gamma\ge \frac{7\pi}{5}$, любого $x\in H$ и любого $\sigma>0$ выполняется неравенство (\ref{vas0}) (с заменой $\widetilde{V}$ на $V^*$) и(\ref{vas}).
 \end{theor}

Теперь рассмотрим вес
\[
\widehat{V}(t)=\left\{ \begin{array}{clcr}\mbox{$5/4$}, &\mbox{$t\in\left[0, 1/2\right]$},\\
\mbox{$3/4$},&\mbox{$t\in\left({1}/{2}, 1\right]$}.
\end{array}\right. \;\;\| \widehat{V}\|_1=1.
\]
В~\cite[теорема 2]{Vas} показано, что для $\psi\in \Psi$ и  $\gamma=\pi$ выполняется неравенство
$$
{\cal G}(\widehat{V},\frac \gamma\sigma,\sigma)\geq
\frac 34{\cal I}\left(\psi\right).
$$
Из (\ref{est}) выводим, что справедлива
\begin{theor}
 Для любой функции $\varphi :\mathbb{C}\to \mathbb{C}$ такой, что $\psi\in \Psi$,  любого $x\in H$ и любого $\sigma>0$ выполняется неравенство
 \begin{equation}\label{vas1}
 E_{\sigma}\left(x\right)\leq \left(\frac 43\right)^{1/2}\frac{1}{\sqrt{{\cal I}\left(\psi\right)}}
  \omega_{\varphi}\left(x; L_{p,\widehat{V}}\left(\left[0,\frac{\pi}{\sigma}\right]\right)\right),\qquad 2\le p\le\infty .
\end{equation}
 \end{theor}

\bigskip

{\bf 7. Некоторые приложения.} Приведенные выше результаты включают в себя ряд точных неравенств типа Джексона -- Стечкина (см., например, \cite{Vas, Vasil}) для наилучших $L_2$ - приближений периодических функций тригонометрическими полиномами, результаты по наилучшим $L_2$ приближениям функций, заданных на всей оси целыми функциями экспоненциального типа, а также аналогичные результаты для почти периодических функций. Уже в этих конкретных случаях некоторые результаты являются новыми.

Результаты для периодических функций получаются, если в пространстве $L_2(\mathbb{T})$ выбрать разложение единицы следующим образом. Для любоых борелевского множества $\beta\subset\mathbb{T}$ и функции $x\in L_2(\mathbb{T})$
\[
E(\beta)x(t):=\sum\limits_{k\in \beta}\hat{x}_ke^{ikt}.
\]

Результаты для функций из $L_2(\mathbb{R})$ получаются, если в этом пространстве выбрать разложение единицы, отвечающее оператору $i\frac d{dt}$, для которого, в силу формулы обращения преобразования Фурье,
\[
E([s,t])x(t)=\frac{1}{2\pi}\int\limits _{-\infty}^{+\infty}\frac{e^{it(z-u)}-e^{is(z-u)}}{i(z-u)}x(z)dz,\; s,t\in\mathbb{R},\:s<t.
\]

Подробнее остановимся на случае почти периодических функций.   В линейном пространстве $\Pi$ функций, почти периодических по Бору, можно ввести скалярное произведение (см. \cite[Дополнение, \S7]{Demidovich})
$
(x,y)=M[x(t)\overline {y(t)}]
$
 и норму $\|\cdot\|=(\cdot,\cdot)^{1/2}$. Пополняя полученное предгильбертово пространство, получаем гильбертово пространство $\widetilde{\Pi}$, которое содержится в пространстве $B^2$ функций почти периодических по Безиковичу (пространства такого типа рассматривались, например, в \cite{Kuzmina}). Для $x\in\widetilde{\Pi}$ и $\lambda\in\mathbb{R}$ положим
\[
a(x,\lambda)=M[x(t)e^{-i\lambda t}].
\]
Как известно, для любой $x\in\widetilde{\Pi}$ множество
$
{\rm Sp}(x)=\{\lambda \; :\; a(x,\lambda)\neq 0\}
$
не более чем счетно.

В пространстве $\widetilde{\Pi}$ определим разложение единицы следующим образом. Для любых борелевского множества $\beta\subset\mathbb{R}$ и функции $x\in \widetilde{\Pi}$
\[
E(\beta)x(t):=\sum\limits_{\lambda\in \beta\cap {\rm Sp}(x)}a(x,\lambda)e^{i\lambda t}.
\]
Тогда будем иметь
\[
\int\limits_{-\infty}^\infty
dE(s)x(t)=\sum\limits_{\lambda\in {\rm Sp}(x)}a(x,\lambda)e^{i\lambda t}.
\]
Для $\sigma>0$ положим
\[
{\cal E}_\sigma(x)=\inf\limits_{y\in\widetilde{\Pi}\atop {\rm Sp}(y)\subset (-\sigma,\sigma)}\| x-y\|.
\]
Будем для простоты рассматривать только функции $x\in \widetilde{\Pi}$, для которых ${\rm Sp}(x)$ имеет единственную предельную точку в бесконечности.
В силу равенства Парсеваля
\[
{\cal E}^2_\sigma(x)=\sum\limits_{\lambda\in {\rm Sp}(x)\atop |\lambda|\ge \sigma}|a(x,\lambda)|^2.
\]
Теперь легко видеть, что для величин ${\cal E}_\sigma(x)$ и  обобщенных модулей непрерывности $\omega_\varphi(x,\delta)$, построенных с помощью определенного выше разложения единицы, справедливы утверждения теоремы 3, следствия 2, а также (при дополнительных предположениях о функции $\varphi$) утверждения теорем 4, 5 и 6.

\end{document}